\documentclass[11pt]{amsart}
\usepackage{amssymb}

\allowdisplaybreaks[1]

\newtheorem{thm}{Theorem}[section]

\newtheorem{lem}[thm]{Lemma}
\newtheorem{prop}[thm]{Proposition}
\newtheorem{defn}[thm]{Definition}
\theoremstyle{definition}

\numberwithin{equation}{section}

\newcommand{\kN}{\mathbb{N^{*}}}
\newcommand{\kkN}{\mathbb{N}}
\newcommand{\kZ}{\mathbb{Z}}
\newcommand{\kA}{\mathbb{A}}

\newcommand{\kC}{\mathbb{C}}
\newcommand{\kI}{\mathbb{I}}
\newcommand{\kR}{\mathbb{R}}
\begin{document}
\title{A class of generalized gamma functions}
\author{Jean-Paul JURZAK}
\vskip.5cm

\begin{abstract}
In this paper, we study the holomorphic function defined by the infinite product $\Gamma_{a,r}(s) =\prod_{ n \geq 0} ( 1 + \frac{1}{a+ nr})^s\;( 1 + \frac{s}{a+nr})^{-1}$ which generalize Euler's definition in the sense that $\Gamma(s) = \frac{\Gamma_{1,1}(s)}{s}$. We obtain analogues of the Gauss multiplication formula, complement formula for functions $\Gamma_{a,r}(s)$.
\end{abstract}
\maketitle
\begin{center}
{  {\small{Laboratoire Gevrey de Math\'ematique Physique, Universit\'e de Bourgogne, \\
Facult\'e des Sciences et Techniques\\
 BP 47870, 21078 Dijon Cedex\\
 e-mail : Jean-Paul.Jurzak@u-bourgogne.fr}}
}
\end{center}

\section{Introduction} \label{sec:intro}
Analogues of the gamma function $\Gamma(s)$ have been proposed and studied by L.Bendersky \cite{Bend33},  E.L.Post\cite{post19}, S.Ramanujan \cite{ramanu} and other authors such as  F.H.Jackson  which introduces a $q$-gamma function defined as an infinite product.
 
Euler's original definition of $\Gamma(s)$ is
\begin{equation} \label{f:Eulerdef}
\Gamma(s) \;=\; \frac{1}{s}\;\prod_{ n \in \kN} \frac{( 1 + \frac{1}{n})^s}{( 1 + \frac{s}{n})} 
\end{equation}
and many representations of $\Gamma(s)$ as integrals involving classical functions can be found in classical manuals (see \cite{tables65} for example).

Our paper differs from previous approachs and study a modified gamma function discribed by formula $\Gamma_{a,r}(s) =\prod_{ n \geq 0} ( 1 + \frac{1}{a+ nr})^s\;( 1 + \frac{s}{a+nr})^{-1}$ closed to Euler's definition. However, no simple integral formula (involving classical functions distinct from Euler gamma function) has been obtained as an alternative description of the function $\Gamma_{a,r}(s)$.

The paper is organized as follows: main properties of $\Gamma_{a,r}(s)$ are summed in proposition \ref{prop:main} and proofs of results of this section are concentrated in section \ref{sec:proofs}. Section \ref{sec:prelimi} states few properties of  infinite products
 defined by formula \ref{ff:gammaAdef}.

\section{Preliminaries} \label{sec:prelimi}
For a real sequence $ \kA=(a_n)_{n \geq 0}$ such that $ 1 + \frac{1}{a_n} >0 $ for all $n$, we formally put
\begin{equation} \label{ff:gammaAdef}
\Gamma_{\kA}(s) \;=\; \prod_{ n \geq 0} \frac{( 1 + \frac{1}{a_n})^s}{( 1 + \frac{s}{a_n})}
\end{equation}
Clearly $\Gamma_{\kA}(0)= \Gamma_{\kA}(1)=1$.
\begin{defn}
Let $ \kA=(a_i)_{i \in \kI}$ be a subset of $\kR$. For $\alpha \in \kC^{*}$ and $n \in \kZ$, we define
\begin{center}
$\frac{\kA}{\alpha}=(\frac{a_i}{\alpha })_{i \in \kI} $ ;
$\kA + \alpha = (a_i+ \alpha)_{i \in \kI}$ ;
$\kA^n = (a_i^{n})_{i \in \kI}$
\end{center}
\end{defn}

\medskip

\begin{prop} \label{prop:mero}
We assume that $\Sigma_{n=1}^{\infty} \frac{1}{a_n^2}\;<\; + \infty$.  One has the following properties:
\begin{enumerate}
\item  The function $ \Gamma_{\kA}(s)$ is defined and meromorphic on $\kC$ with poles only at $s=-a_n$ for $n \in \mathbb{N}$ and never vanishes on $\kC - (-\kA)$.
\item For $\alpha \in \kR^{*}$ and $s \in \kC - (-\frac{\kA}{\alpha})$, one has:
\begin{equation}  \label{f:gamalph}
\Gamma_{\kA}(\alpha\;s) \;=\;\Gamma_{\kA}(\alpha)^s \;\Gamma_{\frac{\kA}{\alpha}}(s)
\end{equation} 
\item Assuming $ a_n \not=  \pm 1$ for all $n$, one has, for $s \in \kC - (\pm\kA)$
\begin{equation}   \label{f:gas2moinss}
 \Gamma_{\kA}(s) = 
 \mu_{\kA}^{s-s^2}\;\frac{ \Gamma_{-\kA}(s)}{ \Gamma_{-\kA^2}(s^2)}
\end{equation}
with $\mu_{\kA} = \prod_{ n \geq 0} ( 1 - \frac{1}{a_n^2})$.
\end{enumerate}
\end{prop}
In this context, the strictly positive real number  
$e^{-\gamma_A}= \prod_{ n \geq 0} (1+\frac{1}{a_n})e^{- \frac{1}{a_n}} $ is  well defined and  one has, for  $s \in \kC - (-\kA)$:
\begin{equation} \label{f:Dgaga}
\frac{\Gamma_A'(s) }{\Gamma_A(s) }
=
-\gamma_A + \Sigma_{n=0}^{\infty}\;\frac{s}{a_n(a_n + s)}
\end{equation}
\medskip
Choosing $\kA=\kN$, we get $ \Gamma(s) \;=\; \frac{1}{s}\;\Gamma_{\kN}(s)$ with $\gamma_A$ 
 equal to the Euler constant $\gamma$.

\textit{Proof:}\\
Putting $J_n = \mathbb{N} \cap [0,n] $ for $n \in \mathbb{N}$, one has, for $s \in \kR$:
\begin{equation} \label{f:gaexpo}
 \prod_{ k \in J_n} \frac{( 1 + \frac{1}{a_k})^s}{( 1 + \frac{s}{a_k})}\;=\;
 \frac{\prod_{ k \in J_n}\big(( 1 + \frac{1}{a_k})\;e^{-\frac{1}{a_k}}\big)^s}
{\prod_{ k \in J_n}( 1 + \frac{s}{a_k})e^{-\frac{s}{a_k}}} 
\end{equation}
The sequence of entire functions $s \rightarrow \prod_{ k \in J_n}( 1 + \frac{s}{a_k})e^{-\frac{s}{a_k}}$ is  uniformly convergent (as $ n \rightarrow +\infty $) on compact sets of $\kC$ to the entire function
$s \rightarrow \prod_{ n \in \kN}( 1 + \frac{s}{a_k})e^{-\frac{s}{a_k}}$ which vanishes only at points $s$ satisfying $( 1 + \frac{s}{a_k})e^{-\frac{s}{a_k}}=0$, ie for $ s \in \kA$. Clearly 
$$ \prod_{ k \in J_n}\big(( 1 + \frac{1}{a_k})\;e^{-\frac{1}{a_k}}\big)^s 
\quad \rightarrow \quad e^{-\gamma_{\kA}s} $$
as $ n$ tends to $+\infty$, uniformly for $s$ in a compact set of $\kC$. Thus, formula 
\ref{ff:gammaAdef} defines a meromorphic function on $\kC$ with poles only at $s=-a_k$ for $n \geq 0$, as a quotient of two entire functions.\\
For formula \ref{f:gamalph}, omitting terms of the form $e^{\pm\frac{1}{a_k}}$, it suffice to prove the equality for $s\in \kR^{+*}$ and formula holds for a general $s \in \kC - (-\frac{\kA}{\alpha})$ by analytic continuation. One has, for $\alpha \in \kR^{+*}$
$$ \frac{ \Gamma_{\kA}(\alpha\;s)}{ \Gamma_{\kA}(\alpha)^s} \;=\;
 \text{lim}_{n \rightarrow \infty}\;
\prod_{ k \in J_n} \frac{( 1 + \frac{1}{a_k})^{\alpha s}}{( 1 + \frac{\alpha s}{a_k})}
\; \bigg[
\prod_{ k \in J_n} \frac{( 1 + \frac{\alpha}{a_k})}{( 1 + \frac{1}{a_k})^{\alpha}}
\bigg]^s $$
$$ = \text{lim}_{n \rightarrow \infty}\;
\prod_{ k \in J_n} \frac{( 1 + \frac{\alpha}{a_k})^{ s}}{( 1 + \frac{\alpha s}{a_k})}
\;=\; \Gamma_{\frac{\kA}{\alpha}}(s)
$$
For formula \ref{f:gas2moinss}, one has:
$$ \Gamma_{\kA}(s) = \prod_{ n \geq 0} \frac{( 1 + \frac{1}{a_n})^s}{( 1 + \frac{s}{a_n})} 
= 
\prod_{ n \geq 0} \frac{( 1 - \frac{1}{a_n^2})^s}{( 1 - \frac{1}{a_n})^s} 
\frac{( 1 - \frac{s}{a_n})}{( 1 - \frac{s^2}{a_n^2})}$$
$$
= 
\prod_{ n \geq 0} \frac{( 1 - \frac{1}{a_n^2})^{s^2}}{( 1 - \frac{1}{a_n})^s} 
\frac{( 1 - \frac{s}{a_n})}{( 1 - \frac{s^2}{a_n^2})}
\frac{( 1 - \frac{1}{a_n^2})^s}{( 1 - \frac{1}{a_n^2})^{s^2}} =
 \mu_{\kA}^{s-s^2}\;\frac{ \Gamma_{-\kA}(s)}{ \Gamma_{-\kA^2}(s^2)}
$$

\vspace{0.6cm}
We note the
\begin{prop}
The functional equation satisfied by the Riemann zeta function is:
\begin{equation} \label{f:zetafonct}
\frac{  \Gamma_{2\kN}(s)\;\zeta(s)} {2^s \;s} \;=\;
\frac{ \Gamma_{2\kN}(1-s)\;\zeta(1-s)} {2^{1-s} \;(1-s)}
\end{equation}
\end{prop}
This equality is the starting point of this work.

\textit{Proof:}\\
Since $\Gamma(\frac{1}{2})=\sqrt{\pi}$, the functional equation takes the form
$$\frac{  \Gamma(\frac{s}{2})\;\zeta(s)}{\Gamma(\frac{1}{2})^s}
\;=\;
\frac{  \Gamma(\frac{1-s}{2})\;\zeta(1-s)}{\Gamma(\frac{1}{2})^{1-s}}$$
hence, from $\Gamma(s) \;=\; \frac{1}{s}\;\Gamma_{\kN}(s)$
$$\frac{  \Gamma_{\kN}(\frac{s}{2})\;\zeta(s)}{2^s\;s\;\Gamma_{\kN}(\frac{1}{2})^s}
\;=\;
\frac{  \Gamma_{\kN}(\frac{1-s}{2})\;\zeta(1-s)}{2^{1-s}\;(1-s)\;\Gamma_{\kN}(\frac{1}{2})^{1-s}}$$
Taking $\kA=\kN$ and $\alpha=\frac{1}{2}$ in formula \ref{f:gamalph} gives
$$ \frac{\Gamma_{\kN}(\frac{s}{2}) }{\Gamma_{\kN}(\frac{1}{2})^s } \;=\; 
\Gamma_{2\kN}(s)  $$
A direct substitution imply formula \ref{f:zetafonct}.

\section{Properties of $\Gamma_{a,r}$ functions}
Let $\kA=(a+nr)_{n \geq 0}= a+r\kkN$ be an arithmetical sequence with reason $r>0$ and first element $a>0$ (or with reason $r<0$ and first element $a<-1$). We define
\begin{equation} \label{f:defgammaA}
\Gamma_{a,r}(s) \;=\; \prod_{ n \geq 0} \frac{( 1 + \frac{1}{a+ nr})^s}{( 1 + \frac{s}{a+nr})}
\end{equation}
This formula is closed to Euler's definition \ref{f:Eulerdef} of gamma function expressed as an infinite product. The $\Gamma_{a,r}$ function is expressed by \ref{f:gammaA} in a non-natural way in terms of the well-known  $\Gamma$ function in order to obtain shorter proofs of some properties of 
$\Gamma_{a,r}$  functions. Some properties of $\Gamma_{a,r}$ are more natural from explicit formula \ref{f:defgammaA}. We note that $ \Gamma(s) \;=\; \frac{1}{s}\;\Gamma_{1,1}(s)$ and 
formula \ref{f:zetafonct} becomes
\begin{equation} \label{f:zeta22fonct}
\frac{  \Gamma_{2,2}(s)\;\zeta(s)} {2^s \;s} \;=\;
\frac{ \Gamma_{2,2}(1-s)\;\zeta(1-s)} {2^{1-s} \;(1-s)}
\end{equation}

Following section \ref{sec:prelimi}, we put
\begin{equation}  \label{defsinar}
\sin_{a,r}(s) \;=\; \prod_{ n \geq 0}( 1 - \frac{s^2}{(a+nr)^2})\;=\;
\frac{1}{\Gamma_{a,r}(s)\Gamma_{a,r}(-s) } 
\end{equation}
$$ e^{-\gamma_{a,r}} = \prod_{ n \geq 0} (1+\frac{1}{a+nr})e^{- \frac{1}{a+nr}}
$$

\medskip

\begin{lem} \label{lem:cstes}
Let $\Psi(s)=\frac{\Gamma'(s) }{\Gamma(s) }$. One has,  for $s \in \kC - \kA$:
\begin{equation} \label{f:sinA}
\sin_{a,r}(s) \;=\; \frac{\Gamma^2(\frac{a}{r})}
{\Gamma(\frac{a}{r}+\frac{s}{r}  )\;\Gamma(\frac{a}{r}-\frac{s}{r}  )}
\end{equation}
\begin{equation}  \label{f:smallgammaA}
e^{-\gamma_{a,r}} = 
\frac{\Gamma(\frac{a}{r})}{\Gamma(\frac{a}{r}+\frac{1}{r})}\;e^{\frac{1}{r} \Psi( \frac{a}{r}) }
\end{equation}
\begin{equation} \label{f:gammaA}
\Gamma_{a,r}(s) = \frac{ 1 }{ \Gamma( \frac{a}{r})} \;\Gamma( \frac{s+a}{r})\;
\frac{ \Gamma( \frac{a}{r})^s}{ \Gamma( \frac{a}{r}+\frac{1}{r})^s } 
\end{equation}
\end{lem}

\medskip

\begin{prop}  \label{prop:main}
The function $\Gamma_{a,r}(s)$ has the following properties 
\begin{enumerate}
\item For $s \in \kC - \kA$
\begin{equation} \label{f:multi}
 \Gamma_{a,r}(s+r) = \alpha^r\;\frac{s+a}{r}\;\Gamma_{a,r}(s)
\end{equation}
with $\alpha=\frac{ \Gamma( \frac{a}{r})}{ \Gamma( \frac{a}{r}+\frac{1}{r}) } $.
\item Complement's formula
\begin{equation} \label{f:compl}
\Gamma_{a,r}(s)\Gamma_{a,r}(r-s)
\;=\;(1-\frac{s}{a})\;
\frac{ \Gamma_{a,r}(r) }{\sin_{a,r}(s) } 
\end{equation}
\item Duplication formula  
\begin{equation}  \label{f:dupli}
(2^{2s})^{ \frac{1}{r}}\;\Gamma_{a,r}( s )\;
\Gamma_{a,r}( s +\frac{r}{2})  \;=\;
\frac{\Gamma_{a,r}( \frac{r}{2}  )}{\Gamma_{a,r}( a )}\;\Gamma_{a,r}( 2s + a)
\end{equation}
with 
$$\frac{\Gamma_{a,r}( \frac{r}{2}  )}{\Gamma_{a,r}( a )}
= \frac{\sqrt{ \pi } }{ 2^{ \frac{2a}{r}-1} }\;
\frac{\Gamma(\frac{a}{r}+\frac{1}{r})^{ a - \frac{r}{2} } }
{\Gamma(\frac{a}{r})^{1+a - \frac{r}{2} } }
$$
\item Multiplication formula
$$  n^{ \frac{ns}{r}}\;
\prod_{k=0}^{k=n-1}\; \Gamma_{a,r}(s+k\frac{r}{n})=
\frac{\Gamma_{a,r}( \frac{r}{n}) \cdots \Gamma_{a,r}( \frac{kr}{n}) \cdots \Gamma_{a,r}( \frac{(n-1)r}{n}) }{\Gamma_{a,r}( (n-1)a )}
\;\Gamma_{a,r}(ns+(n-1)a)
 $$
with
\begin{equation} \label{f:ctesmultip}
\frac{\Gamma_{a,r}( \frac{r}{n}) \cdots \Gamma_{a,r}( \frac{kr}{n}) \cdots \Gamma_{a,r}( \frac{(n-1)r}{n}) }{\Gamma_{a,r}( (n-1)a )}
=\frac{ (2\pi)^{ \frac{n-1  }{ 2 } } }
{ n^{ \frac{na}{r}-\frac{1}{2}} }
\;
\frac{\Gamma(\frac{a}{r}+\frac{1}{r})^{ (n-1)(a - \frac{r}{2})} }
{\Gamma(\frac{a}{r})^{(n-1)(a +1 - \frac{r}{2}) } }
 \end{equation}
\end{enumerate}
\end{prop}

Thus, the triple formula may be written as 
$$\Gamma_{a,r}(3s+2a)= \text{Cte} \; 3^{ \frac{3s}{r}}\;
\Gamma_{a,r}(s)\;\Gamma_{a,r}(s+\frac{r}{3})\;\Gamma_{a,r}(s+\frac{2r}{3})$$
with
$$ \text{Cte} =\frac{ 3^{ \frac{3a}{r}-\frac{1}{2}} }{\sqrt{ 2\pi } }\;
\frac{\Gamma(\frac{a}{r})^{2+2a - r } }{
\Gamma(\frac{a}{r}+\frac{1}{r})^{ 2a - r} }
=\frac{\Gamma_{a,r}(2a) }{\Gamma_{a,r}(\frac{r}{3}) \Gamma_{a,r}(\frac{2r}{3})}
 $$
We put $\Psi_{a,r}(s) = \frac{\Gamma_{a,r}'(s) }{\Gamma_{a,r}(s) }$. One has,  for $s \in \kC - (a + r\kkN) $:
$$ \Psi_{a,r}(s) =  \frac{1}{r}\; \Psi(\frac{s+a}{r}) +
\ln\Gamma( \frac{a}{r}) - \ln\Gamma( \frac{a}{r}+\frac{1}{r})$$
$$ \Psi_{a,r}(s+r) -  \Psi_{a,r}(s) = \frac{1}{s+a}$$

\vspace{0.6cm}
Various properties are
\begin{prop}
Let $\beta(\cdot,\cdot)$ be the usual beta function. We assume $a >0$ and $r>0$. With obvious notations, one has:
\begin{enumerate} \label{propvarious}
\item For $ n \in \kN$
\begin{equation} \label{ggwithnr}
\Gamma_{a,r}(s) = \Gamma_{a,nr}(s)\;\Gamma_{a+r,nr}(s)\;\cdots\Gamma_{a+kr,nr}(s)\;\cdots \Gamma_{a+(n-1)r,nr}(s) 
\end{equation}
and
$$ \frac{  \Gamma_{a,nr}(s) }{ \Gamma_{a+kr,nr}(s)  }  = 
 \frac{\beta( \frac{s+a}{nr} ; \frac{k}{n} )}{\beta( \frac{a}{nr} ; \frac{k}{n} ) }\;
 \frac{\beta( \frac{a}{nr} ; \frac{k}{n} )^s}{\beta( \frac{a+1}{nr} ; \frac{k}{n} )^s }
=\frac{\beta( \frac{s}{nr} ; \frac{a+kr}{nr} )}{\beta( \frac{s}{nr} ; \frac{a}{nr} ) }\;
 \frac{\beta( \frac{a}{nr} ; \frac{1}{nr} )^s}{\beta( \frac{a+kr}{nr} ; \frac{1}{nr} )^s }
$$
\item For suitable $s$ and $h$, one has
\begin{equation} \label{f:amoinsh} 
 \Gamma_{a-h,r}(s+h) =  C \; q^s\; \Gamma_{a,r}(s)  
\end{equation}
with
$$ C = \Gamma(\frac{a}{r})\;
  \frac{\Gamma(\frac{a-h}{r})^{h-1} }{\Gamma(\frac{a-h}{r}+\frac{1}{r})^h } 
\quad \text{and} \quad 
 q = \frac{\beta( \frac{a-h}{r} ; \frac{1}{r} )}{\beta( \frac{a}{r} ; \frac{1}{r} ) }
=\frac{\beta( \frac{a+1}{r} ; \frac{-h}{r} )}{\beta( \frac{a}{r} ; \frac{-h}{r} ) }
$$
\item One has
\begin{equation} \label{f:newgamm}
\Gamma_{a,r}(s) \;= \;\frac{1}{\Gamma( \frac{a}{r}) }\;\frac{\Gamma_{a,r}(r)^\frac{s}{r} }{(\frac{a}{r})^\frac{s}{r} } \;\Gamma( \frac{s}{r}+\frac{a}{r})
\end{equation}
\end{enumerate}
\end{prop}

\vspace{0.6cm}
When $a=r$, the function $\Gamma_{a,a}(s)$  becomes
\begin{equation} \label{f:defgamma}
\Gamma_{a,a}(s) \;=\;\prod_{ n \in \kN} \frac{( 1 + \frac{1}{an})^s}{( 1 + \frac{s}{an})} 
\end{equation}
with properties
\begin{equation} \label{f:asinA}
\sin_{a,a}(s) \;=\; a\;\frac{\sin(\pi\frac{s}{a})}
{\pi\; s}
\end{equation}
\begin{equation}  \label{f:asmallgammaA}
e^{-\gamma_{a,a}} =_{def} 
\prod_{ n \in \kN} (1+\frac{1}{an})e^{- \frac{1}{an}}
=
\frac{e^{-\frac{\gamma}{a}}}
{\Gamma(1+\frac{1}{a})} \qquad \gamma \quad \text{Euler constant}
\end{equation}
\begin{equation} \label{f:agammaA}
\Gamma_{a,a}(s) = 
\frac{ \Gamma( 1+\frac{s}{a})}{ \Gamma( 1 +\frac{1}{a})^s } 
\end{equation}
\begin{equation} \label{f:amulti} 
 \Gamma_{a,a}(s+a) = \frac{1}{ \Gamma(1+\frac{1}{a})^a }\;\frac{s+a}{a}\;\Gamma_{a,a}(s)
\end{equation}
\begin{equation} \label{f:acompl}
\Gamma_{a,a}(s)\Gamma_{a,a}(a-s)
\;=\;\frac{1}{ \Gamma(1+\frac{1}{a})^a }(1-\frac{s}{a})\;
\frac{\pi\frac{s}{a}}{\sin(\pi\frac{s}{a})}
\end{equation}
\begin{equation} \label{f:gagagbgb}
\Gamma_{a,a}(\frac{a}{b}\;s) \;=\;\Gamma_{a,a}(\frac{a}{b})^s \;\Gamma_{b,b}(s) 
\end{equation}
\begin{equation} \label{f:aapsi}
 \Psi_{a,a}(s) =  \frac{1}{a}\; \Psi(1+\frac{s}{a})  - \ln\Gamma( 1+\frac{1}{a})
\end{equation}

\medskip

\begin{prop} \label{propsin}

For $\lambda \neq 0$
\begin{equation} \label{sinlambda}
\sin_{a,r}(s) \;=\;\sin_{\lambda a,\lambda r}(\lambda s)
\;=\;\sin_{\frac{a}{r},1}(\frac{ s}{r})
\end{equation}
\begin{equation}  \label{sinsplr}
\frac{ \sin_{a,r}(s+r)}{(\frac{s}{r} - \frac{a}{r}+1) } \;=\;-\;\frac{\sin_{a,r}(s)}{ (\frac{s}{r}+\frac{a}{r}  )}
\end{equation}
\begin{equation}   \label{sinspr2}
\frac{ \sin_{a,r}(s+\frac{r}{2})}{(\frac{s}{r}-\frac{a}{r}+\frac{1}{2}  ) } \;=\;-\; 
\frac{\Gamma^2(\frac{a}{r})  }
{\Gamma^2(\frac{a}{r}+\frac{1}{2}) }\;\sin_{a+\frac{r}{2},r}(s) 
\end{equation}
\begin{equation}  \label{sinamoinsr2}
 \sin_{a-\frac{r}{2},r}(s) \;=\;(1 - \frac{s^2}{(a-\frac{r}{2} )^2} )\;
\sin_{a+\frac{r}{2},r}(s)
\end{equation}
For $ n \in \kN$
\begin{equation} \label{sinwithnr}
\sin_{a,r}(ns) = \sin_{\frac{a}{n},r}(s)\;
\sin_{\frac{a}{n}+\frac{r}{n},r}(s)\;\cdots
\sin_{\frac{a}{n}+\frac{kr}{n},r}(s)\;\cdots
\sin_{\frac{a}{n}+\frac{(n-1)r}{n},r}(s)
\end{equation}
\end{prop}

\vskip18pt

Taking $n=2$ in formula \ref{sinwithnr}, we find
\begin{equation} \label{sintwos}
\sin_{a,r}(2 s) \;=\; \sin_{\frac{a}{2},r}(s) \;
\sin_{\frac{a}{2}+\frac{r}{2},r}(s)
\end{equation}
Classical descriptions of $\sin(s)$ and $\cos(s)$ as infinite product are equivalent to
$$ \frac{sin(\pi \;s)}{\pi \;s}=sin_{1,1}(s) \qquad \text{et} \qquad cos(\pi\;s)=sin_{\frac{1}{2},1}(s)$$
with the following properties
$$ sin_{\frac{1}{2},1}(s)=\pi(s + \frac{1}{2})\;sin_{1,1}(s+ \frac{1}{2})$$
$$ sin_{\frac{1}{2},1}(s+1) = -sin_{\frac{1}{2},1}(s)$$
due to $cos(\pi\;s)=\sin(\pi\;s + \frac{\pi}{2})$ and periodicity of $\cos(s)$. And $sin(2\pi s)=2\;sin(\pi s)\;cos(\pi s)$ expresses as 
$$sin_{1,1}(2s)=sin_{1,1}(s)\;sin_{\frac{1}{2},1}(s)= 
\pi(s + \frac{1}{2})\;sin_{1,1}(s)\;sin_{1,1}(s+ \frac{1}{2})
$$
which agree with formulas \ref{sintwos} and \ref{sinspr2} choosing $a=1$ and $r=1$.

\medskip

\section{Proofs}  \label{sec:proofs}
\textit{Proof of lemma \ref{lem:cstes}}:\\
One has from \cite{tables65} formula 8.364 
$$\prod_{ n \geq 0}\; (1+ \frac{w}{a+nr})\;e^{- \frac{w}{a+nr}}
= \frac{ \Gamma( \frac{a}{r})}{ \Gamma( \frac{a}{r}+\frac{w}{r}) } \; e^{\frac{w}{r}\; \Psi (\frac{a}{r})}$$
hence $w=1$ gives formula \ref{f:smallgammaA}. Choosing $w=\pm s$ gives formula \ref{f:sinA}. For $s \in \kR$
$$\Gamma_{a,r}(s) = \prod_{ n \geq 0}\; \frac{\bigg[(1+ \frac{1}{a+nr})\;e^{- \frac{1}{a+nr}}\bigg]^s }
{(1+ \frac{s}{a+nr})\;e^{- \frac{s}{a+nr}}}
= \frac{ \frac{ \Gamma( \frac{a}{r})^s}{ \Gamma( \frac{a}{r}+\frac{1}{r})^s } \; e^{\frac{s}{r}\; \Psi (\frac{a}{r})}
 } {\frac{ \Gamma( \frac{a}{r})}{ \Gamma( \frac{a}{r}+\frac{s}{r}) } \; e^{\frac{s}{r}\; \Psi (\frac{a}{r})} }
$$
thus, by analytic continuation and proposition  \ref{prop:mero}
$$\Gamma_{a,r}(s) = \frac{ \frac{ \Gamma( \frac{a}{r})^s}{ \Gamma( \frac{a}{r}+\frac{1}{r})^s } \; 
 } {\frac{ \Gamma( \frac{a}{r})}{ \Gamma( \frac{a}{r}+\frac{s}{r}) }  } =
\frac{ \Gamma( \frac{a}{r}+\frac{s}{r}) }{ \Gamma( \frac{a}{r})} \;
\frac{ \Gamma( \frac{a}{r})^s}{ \Gamma( \frac{a}{r}+\frac{1}{r})^s }
$$

\medskip

\textit{Proof of proposition  \ref{prop:main}}:\\
Let $\alpha=\frac{ \Gamma( \frac{a}{r})}{ \Gamma( \frac{a}{r}+\frac{1}{r}) } $. One has:
$$\frac{\Gamma_{a,r}(s)  }{\Gamma_{a,r}(s+r)   }= \frac{1}{\alpha^r}\; \frac{ \Gamma(\frac{s+a}{r} ) }
{ \Gamma(\frac{s+r+a}{r} ) }
=\frac{1}{\alpha^r}\; \frac{ \Gamma(\frac{s+a}{r} ) }
{ \frac{s+a }{r }\Gamma(\frac{s+a}{r} ) }
=\frac{1}{\alpha^r}\;\frac{r}{s+a}$$
For complement formula \ref{f:compl}, with $J_n = \mathbb{N} \cap [0,n] $ for $n \in \mathbb{N}$, one has, for $s \in \kR$:
$$ \Gamma_{a,r}(s)\Gamma_{a,r}(r-s)
 \;=\;  \text{lim}_{n \rightarrow \infty}\;
\prod_{ k \in J_n} 
\frac{1}{1+\frac{s}{a+kr}}
\frac{(1+\frac{1}{a+kr})^{r}}{(1+\frac{r-s}{a+kr})}
$$
At fixed $n$, the denominator is
$$ \bigg((1+\frac{s}{a})
\frac{(a+r+s)}{(a+r)}\frac{(a+2r+s)}{(a+2r)}
\cdots
\frac{(a+nr+s)}{(a+nr)}\bigg)
\bigg(
\frac{(a+r-s)}{a}\frac{(a+2r-s)}{(a+r)}
\cdots
\frac{(a+(n+1)r-s)}{(a+nr)}\bigg)
$$
$$=\; (1+\frac{s}{a})\bigg[\bigg(\frac{ (a+r)^2-s^2}{(a+r)^2}\bigg)
\bigg(\frac{ (a+2r)^2-s^2}{(a+2r)^2}\bigg)
\cdots
\bigg(\frac{ (a+nr)^2-s^2}{(a+nr)^2}\bigg)
\bigg]\frac{(a+(n+1)r-s)}{a}$$
$$
=\; \frac{1}{(1-\frac{s}{a})}
\bigg[( 1-\frac{s^2}{a^2})
( 1-\frac{s^2}{(a+r)^2})
( 1-\frac{s^2}{(a+2r)^2})
\cdots
( 1-\frac{s^2}{(a+nr)^2}) 
\bigg]\frac{(a+(n+1)r-s)}{a}
$$
Thus
$$\Gamma_{a,r}(s)\Gamma_{a,r}(r-s)\;=\;
 \text{lim}_{n \rightarrow \infty}\;
(1-\frac{s}{a})\; \frac{a }{(a+(n+1)r -s ) }\;
\prod_{ k \in J_n} 
\frac{ (1+\frac{1}{a+kr})^{r} }
{ ( 1-\frac{s^2}{(a+kr)^2}) }
$$
One has 
$$ \prod_{ k \in J_n}  (1+\frac{1}{a+kr})^{r}=
\prod_{ k \in J_n}  \bigg((1+\frac{1}{a+kr})e^{-\frac{1}{a+kr} } \bigg)^{r}
\;
e^{\;\sum_{ k \in J_n} \frac{r}{a+kr}}
$$
$$ \sum_{ k \in J_n} \frac{r}{a+kr}
=
\Psi(1+n + \frac{a}{r}) -  \Psi(\frac{a}{r})
$$
Clearly, as $n \rightarrow +\infty$    
$$ \prod_{ k \in J_n}( 1-\frac{s^2}{(a+kr)^2}) \quad \rightarrow \quad \sin_{a,r}(s)$$
$$ \prod_{ k \in J_n}  \bigg((1+\frac{1}{a+kr})e^{-\frac{1}{a+kr} } \bigg)^{r}
\quad \rightarrow \quad e^{-r\;\gamma_{a,r}}
$$
$$
\Psi(1+n + \frac{a}{r}) = \ln(1+n + \frac{a}{r}) - \frac{1}{2(1+n + \frac{a}{r}) } + O\big(\frac{1}{n^2}\big) 
$$
hence
$$\Gamma_{a,r}(s)\Gamma_{a,r}(r-s)\;=\;(1-\frac{s}{a})\;
\frac{ e^{-r\;\gamma_{a,r}}
 }{\sin_{a,r}(s) }
e^{-\Psi(\frac{a}{r}) }\;
 \text{lim}_{n \rightarrow \infty}\;\prod_{ k \in J_n}
 \frac{1}{ (1+\frac{(n+1)r-s}{a})\;e^{-\; \Psi(1+n + \frac{a}{r})}}  
$$
But
$$ \prod_{ k \in J_n} (1+\frac{(n+1)r-s}{a})\;e^{-\; \Psi(1+n + \frac{a}{r})}=
(1+\frac{(n+1)r-s}{a})\;e^{-\; \ln(1+n + \frac{a}{r})+ 
\frac{1}{2n } + O\big(\frac{1}{n^2}\big)}\quad
\rightarrow \quad \frac{r}{a} 
$$
thus
$$\Gamma_{a,r}(s)\Gamma_{a,r}(r-s)\;=\;\frac{r}{a}(1-\frac{s}{a})\;
\frac{ e^{-r\;\gamma_{a,r}}
 }{\sin_{a,r}(s) }
e^{-\Psi(\frac{a}{r}) }$$
and $s=0$ gives
$$ \Gamma_{a,r}(r)= \frac{r}{a} e^{-r\;\gamma_{a,r}} e^{-\Psi(\frac{a}{r}) } 
$$
showing \ref{f:compl} for $s \in \kR$, hence for $s \in \kC-(a+r\kkN)$ by analytic continuation.\\
We prove multiplication formula. Putting $ x= \frac{s+a}{r}$ and using formula \ref{f:gammaA}, one has
$$
\prod_{k=0}^{k=n-1}\; \Gamma_{a,r}(s+k\frac{r}{n})=
\prod_{k=0}^{k=n-1}\; 
\frac{ 1 }{ \Gamma( \frac{a}{r})} \;\Gamma( \frac{s+a}{r} +\frac{k}{n} )\;
\frac{ \Gamma( \frac{a}{r})^{s +\frac{kr}{n}}}{ \Gamma( \frac{a}{r}+\frac{1}{r})^{s +\frac{kr}{n}}}
$$
$$ =
\frac{1}{\Gamma(\frac{a}{r})^n }\;
\frac{\Gamma(x)\Gamma(x+\frac{1}{n})\cdots\Gamma(x+\frac{k}{n})\cdots\Gamma(x+\frac{n-1}{n})\;\Gamma(\frac{a}{r})^{ns+ \frac{ (n-1)r }{ 2 } }  }
{ \Gamma( \frac{a}{r}+\frac{1}{r})^{ns+ \frac{ (n-1)r }{ 2 } } }
$$
From Gauss formula
$$\Gamma(nx)
= (2\pi)^{\frac{1-n}{2} }\; n^{nx- \frac{1}{2} }\;\prod_{k=0}^{k=n-1}\; \Gamma( x+\frac{k}{n}) 
$$
and
$$
\Gamma_{a,r}(ns+(n-1)a)\;\frac
{\Gamma(\frac{a}{r}+\frac{1}{r})
^{(a+s)n-a } }{\Gamma(\frac{a}{r} )^{(a+s)n-a-1 } }  = \Gamma(nx)
$$
we get that
$$\prod_{k=0}^{k=n-1}\; \Gamma_{a,r}(s+k\frac{r}{n})
=
\frac{\Gamma(\frac{a}{r})^{an-n+a+1+ \frac{ (n-1)r }{ 2 } }  }
{ (2\pi)^{\frac{1-n}{2} }\; n^{nx- \frac{1}{2} } \; \Gamma(\frac{a}{r}+\frac{1}{r})^{a-na+ \frac{ (n-1)r }{ 2 } } }
\;\Gamma_{a,r}(ns+(n-1)a)
$$
thus
$$(2\pi)^{\frac{1-n}{2} }\; n^{\frac{ns}{r}}\;n^{\frac{na}{r}- \frac{1}{2} } \;
\prod_{k=0}^{k=n-1}\; \Gamma_{a,r}(s+k\frac{r}{n})=
\frac{  \Gamma(\frac{a}{r}+\frac{1}{r})^{a-na+ \frac{ (n-1)r }{ 2 } } }
{\Gamma(\frac{a}{r})^{an-n+a+1+ \frac{ (n-1)r }{ 2 } }  }
\;\Gamma_{a,r}(ns+(n-1)a)
$$
which agree with \ref{f:ctesmultip}. Taking $n=2$, we get formula  \ref{f:dupli}.\\

\medskip
%%\section{Various properties}
\textit{Proof of proposition \ref{propvarious}}:\\
We prove 1) of proposition \ref{propvarious}. Formula \ref{ggwithnr} follows from
formula \ref{f:defgammaA}. From formula \ref{f:gammaA}, one has
$$ \frac{  \Gamma_{a,nr}(s) }{ \Gamma_{a+kr,nr}(s)  }  =
 \frac{ \Gamma( \frac{a+kr}{nr}) }{ \Gamma( \frac{a}{nr})} 
\frac{ \Gamma( \frac{s+a}{nr}) }{\Gamma( \frac{s+a+kr}{nr})  }
\;\;
 \frac{ \Gamma( \frac{a}{nr})^s }{ \Gamma( \frac{a+kr}{nr})^s } 
\frac{ \Gamma( \frac{a+kr+1}{nr})^s }{ \Gamma( \frac{a+1}{nr})^s }$$
Multiplying by $ \Gamma( \frac{k}{n})$ or by $ \Gamma( \frac{s}{nr})$, one has:
$$ \frac{ \Gamma( \frac{a+kr}{nr}) }{ \Gamma( \frac{a}{nr})} 
\frac{ \Gamma( \frac{s+a}{nr}) }{\Gamma( \frac{s+a+kr}{nr})  }=
\frac{\beta( \frac{s+a}{nr} ; \frac{k}{n} )}{\beta( \frac{a}{nr} ; \frac{k}{n} ) }=
\frac{\beta( \frac{s}{nr} ; \frac{a+kr}{nr} )}{\beta( \frac{s}{nr} ; \frac{a}{nr} ) }
$$
In the same way
$$\frac{ \Gamma( \frac{a}{nr})^s }{ \Gamma( \frac{a+kr}{nr})^s } 
\frac{ \Gamma( \frac{a+kr+1}{nr})^s }{ \Gamma( \frac{a+1}{nr})^s }=
\frac{\beta( \frac{a}{nr} ; \frac{k}{n} )^s}{\beta( \frac{a+1}{nr} ; \frac{k}{n} )^s }=
\frac{\beta( \frac{a}{nr} ; \frac{1}{nr} )^s}{\beta( \frac{a+kr}{nr} ; \frac{1}{nr} )^s }
$$
For 2), one has
$$  \frac{\Gamma_{a-h,r}(s+h) }{\Gamma_{a,r}(s) } 
=  C\; \bigg[
\frac{\Gamma(\frac{a}{r}+\frac{1}{r}) }{\Gamma(\frac{a}{r}) } \;  \frac{\Gamma(\frac{a-h}{r}) }{\Gamma(\frac{a-h}{r}+\frac{1}{r}) } 
\bigg]^s$$
with
$$ C = \Gamma(\frac{a}{r})\;
  \frac{\Gamma(\frac{a-h}{r})^{h-1} }{\Gamma(\frac{a-h}{r}+\frac{1}{r})^h } $$
With a similar method, we get formula \ref{f:amoinsh}.\\
For 2), let $\kA$ be the sequence $\{a+nr\}_{n \geq 0}$. From \ref{f:gamalph} with $\alpha=r$, we obtain
$$\Gamma_{a,r}(r\;s) \;=\;\Gamma_{a,r}(r)^s \;\Gamma_{\frac{a}{r},1}(s)$$
From formula \ref{f:gammaA}
$$
\Gamma_{\frac{a}{r},1}(s) = \frac{ 1 }{ \Gamma( \frac{a}{r})} \;\Gamma( s+\frac{a}{r})\;
\frac{ \Gamma( \frac{a}{r})^s}{ \Gamma( \frac{a}{r}+1)^s }\;=\;
\frac{ 1 }{(\frac{a}{r})^s \Gamma( \frac{a}{r})} \;\Gamma( s+\frac{a}{r})\; 
$$
thus
$$
\Gamma_{a,r}(r\;s) \;= \;\frac{\Gamma_{a,r}(r)^s }{\Gamma( \frac{a}{r})\;(\frac{a}{r})^s } \;\Gamma( s+\frac{a}{r})
$$
which is formula \ref{f:newgamm}.
\medskip

\textit{Proof of proposition \ref{propsin}:}

One clearly has, for $\lambda\in \kR^{*}$, from formula \ref{f:sinA}
\begin{equation} 
\sin_{a,r}(s) \;=\; \frac{\Gamma^2(\frac{a}{r})}
{\Gamma(\frac{a}{r}+\frac{s}{r}  )\;\Gamma(\frac{a}{r}-\frac{s}{r}  )}
\;=\; \frac{\Gamma^2(\frac{ \lambda a}{ \lambda r})}
{\Gamma(\frac{\lambda a}{\lambda r}+\frac{\lambda s}{\lambda r}  )\;\Gamma(\frac{\lambda a}{\lambda r}-\frac{\lambda s}{\lambda r}  )}
\;=\; \sin_{\lambda a,\lambda r}(\lambda s)
\end{equation}
Taking $\lambda= \frac{1}{r}$, we get formula \ref{sinlambda}.\\
For formula \ref{sinsplr}, one has
$$
\sin_{a,r}(s+r) \;=\; \frac{\Gamma^2(\frac{a}{r})}
{\Gamma(\frac{a}{r}+\frac{s}{r} +1 )\;\Gamma(\frac{a}{r}-\frac{s}{r}-1  )}
 \;=\; \frac{(\frac{a}{r}-\frac{s}{r} -1 )}{ (\frac{a}{r}+\frac{s}{r}  )}\;\frac{\Gamma^2(\frac{a}{r})}
{\; \Gamma(\frac{a}{r}+\frac{s}{r}  )\;\Gamma(\frac{a}{r}-\frac{s}{r}  )}
$$
hence
$$ \sin_{a,r}(s+r) \;=\;-\;\frac{(\frac{s}{r} - \frac{a}{r}+1)}{ (\frac{s}{r}+\frac{a}{r}  )}\;
\sin_{a,r}(s)
$$
From formula \ref{f:sinA}, we obtain
$$ 
\sin_{a+\frac{r}{2},r}(s) \;=\; \frac{\Gamma^2(\frac{a}{r}+\frac{1}{2})}
{\Gamma(\frac{a}{r}+\frac{s}{r} +\frac{1}{2} )\;\Gamma(\frac{a}{r}-\frac{s}{r}+\frac{1}{2}  )}
$$
thus
$$
\sin_{a,r}(s+\frac{r}{2}) \;=\; \frac{\Gamma^2(\frac{a}{r})}
{\Gamma(\frac{a}{r}+\frac{s}{r}+\frac{1}{2}  )\;\Gamma(\frac{a}{r}-\frac{s}{r}-\frac{1}{2}  )}
\;=\; \frac{(\frac{a}{r}-\frac{s}{r}-\frac{1}{2}  )\;\Gamma^2(\frac{a}{r}) }
{\Gamma(\frac{a}{r}+\frac{s}{r}+\frac{1}{2}  )\;\Gamma(\frac{a}{r}-\frac{s}{r}+\frac{1}{2} )}
$$
proving so formula \ref{sinspr2}.\\
Formula \ref{sinamoinsr2} follows directly from definition \ref{defsinar}.\\
Using formula \ref{ggwithnr}, we obtain
\begin{equation} 
\sin_{a,r}(s) = \sin_{a,nr}(s)\;\sin_{a+r,nr}(s)\;\cdots\sin_{a+kr,nr}(s)\;\cdots \sin_{a+(n-1)r,nr}(s) 
\end{equation}
Replacing $s$ by $ns$ and applying formula \ref{sinlambda} with $\lambda = \frac{1}{n}$, we obtain formula \ref{sinwithnr}.

\end{document}